\documentclass{article}
\usepackage{amsfonts}
\usepackage{amssymb}
\usepackage{amsmath}
\usepackage[latin5]{inputenc}
\usepackage[mathscr]{eucal}
\usepackage{eufrak}

\numberwithin{equation}{section}
\input{tcilatex}
\begin{document}

\title{ {\Large {\textbf{ON DUAL TIMELIKE MANNHEIM PARTNER CURVES IN }}}$D$%
{\Large {\textbf{$_{1}^{3}$}}}}
\author{Özcan BEKTAŞ $^{*}$ \and Süleyman ŞENYURT \thanks{%
Ordu University, Faculty of Art and Science, Department of Mathematics,
52750, Perşembe, Ordu, Turkey, ozcanbektas1986@hotmail.com,
senyurtsuleyman@hotmail.com.} }
\date{}
\maketitle

\begin{abstract}
The first aim of this paper is to define the dual timelike Mannheim partner
curves in Dual Lorentzian Space $D_{1}^{3}$, the second aim of this paper is
to obtain the relationships between the curvatures and the torsions of the
dual timelike Mannheim partner curves with respect to each other and the
final aim of this paper is to get the necessary and sufficient conditions
for the dual timelike Mannheim partner curves in $D_{1}^{3}$ .
\end{abstract}

\begin{center}
{\footnotesize \textbf{2000 AMS Subject Classification:}
53B30,51M30,53A35,53A04 }

{\footnotesize \textbf{Keywords:} Mannheim curves, dual Lorentzian Space,
curvature, torsion. }
\end{center}

\section{INTRODUCTION}

\indent In the differential geometry, special curves have an important role.
Especially, the partner curves, i.e., the curves which are related each
other at the corresponding points, have drawn attention of many
mathematicians so far. The well-known of the partner curves is Bertrand
curves which are defined by the property that at the corresponding points of
two space curves principal normal vectors are common. Bertrand partner
curves have been studied in ref. [2,3,5,7,17,22] Ravani and Ku have
transported the notion of bertrand curves to the ruled surfaces and called
Bertrand offsets [16]. Recently, Liu and Wang have defined a new curve pair
for space curves. They called these new curves as Mannheim partner curves:
Let $\alpha $ and $\beta ~$be two curves in th three dimensional Euclidean
space. If there exists a correspondence between the space curves $\alpha $
and $\beta $ such that, at the corresponding points of the curves, the
principal normal lines of $\alpha $ coincides with the binormal lines of $%
\beta $, then $\alpha $ is called a Mannheim curve, and $\beta $ is called a
Mannheim partner curve of $\alpha $. The pair $\left\{ \alpha ,\beta
\right\} $ is said to be a Mannheim pair. They showed that the curve $\alpha
(s)~$is the Mannheim partner curve of $\beta (s^{\ast })~$if and only if the
curvature $k_{1}~$and the torsion $k_{2}~$of $\beta (s^{\ast })$ satisfy the
following equation

\bigskip 
\begin{equation*}
k_{2}^{^{\prime }}=\frac{dk_{2}}{ds^{\ast }}=\frac{k_{1}}{\lambda }%
(1+\lambda ^{2}k_{2}^{\ast })
\end{equation*}

\noindent for some non-zero constants $\lambda .$They also studied the
Mannheim partner curves in the Minkowski 3- space and obtained the necessary
and sufficient conditions for the Mannheim partner curves in $E_{1}^{3}$ [
See 8 and 22 for details]. Moreover, Oztekin and Ergut [15] studied the null
Mannheim curves in the same space. Orbay and Kasap gave [13] new
characterizations of Mannheim partner curves in Euclidean 3-space. They also
studied [12] the Mannheim offsets of ruled surfaces in Euclidean 3- space.
The corresponding characterizations of Mannheim offsets of timelike and
spacelike ruled surfaces have been given \ by Onder and et al [9,10]. New
characterizations of Mannheim partner curves are given in Minkowski 3- space
by Kahraman and et al [6].

In this paper, we study the dual timelike Mannheim partner curves in dual
Lorentzian space $D_{1}^{3}$. Furthermore, we show that the Manheim theorem
is not valid for Mannheim partner curves in $D_{1}^{3}.$Moreover, we give
some new characterizations of the Mannheim partner curves by considering the
spherical indicatrix of some Frenet vectors of the curves.

\noindent

\section{PRELIMINARY}

\indent By a dual number $A$, we mean an ordered pair of the form $\left(
a,a^{\ast }\right) $ for all $a,a^{\ast }\in 
\mathbb{R}
$. Let the set $%
\mathbb{R}
\times 
\mathbb{R}
$ be denoted as $D$. Two inner operations and an equality on $ID=\left\{
\left( a,a^{\ast }\right) \left\vert a,a^{\ast }\in 
\mathbb{R}
\right. \right\} $ are defined as follows:

\indent$\left( i\right) \oplus :D\times D\rightarrow D$, $A\oplus B=\left(
a,a^{\ast }\right) \oplus \left( b,b^{\ast }\right) =\left( a+b,a^{\ast
}+b^{\ast }\right) $ is called the addition in $D$,\newline
\indent$\left( ii\right) \odot :D\times D\rightarrow D$. $A\odot B=\left(
a,a^{\ast }\right) \odot \left( b,b^{\ast }\right) =\left( ab,ab^{\ast
}+a^{\ast }b\right) $is called the multiplication in $D$,\newline
\indent$\left( iii\right) $ $A=B$ iff $a=b$, $a^{\ast }=b^{\ast }$.

\indent If the operations of addition, multiplication and equality on $D=%
\mathbb{R}
\times 
\mathbb{R}
$ with set of real numbers $%
\mathbb{R}
$are defined as above, the set $D$ is called the dual numbers system and the
element $(a,a^{\ast })$ of $D$ is called a dual number. In a dual number $%
A=(a,a^{\ast })\in D$, the real number $a$ is called the real part of $A$
and the real number $a^{\ast }$ is called the dual part of $A$ The dual
number $1=(1,0)$ is called the unit element of multiplication operation $D$
with respect to multiplication and denoted by $\varepsilon $. In accordance
with the definition of the operation of multiplication, it can be easily
seen that $\varepsilon ^{2}=0$. Also, the dual number $A=(a,a^{\ast })\in D$
can be written as $A=a+\varepsilon a^{\ast }$. \newline
\indent The set $D=\{A=a+\varepsilon ^{\ast }a|a,a^{\ast }\in 
\mathbb{R}
\}$ of dual numbers is a commutative ring according to the operations,%
\newline
\indent i) $(a+\varepsilon a^{\ast })+(b+\varepsilon b^{\ast
})=(a+b)+\varepsilon (a^{\ast }+b^{\ast })$\newline
\indent ii)$(a+\varepsilon a^{\ast })(b+\varepsilon b^{\ast
})=ab+\varepsilon (ab^{\ast }+ba^{\ast })$. \newline
\indent The dual number $A=a+\varepsilon a^{\ast }$ divided by the dual
number $B=b+\varepsilon b^{\ast }$ provided $b\neq 0$ can be defined as 
\newline
\indent$\frac{A}{B}=\frac{a+\varepsilon a^{\ast }}{b+\varepsilon b^{\ast }}=%
\frac{a}{b}+\varepsilon \frac{a^{\ast }b-ab^{\ast }}{b^{2}}.$ \newline
\indent Now let us consider the differentiable dual function. If the dual
function $f$ expansions the Taylor series then we have \newline
\indent$f(a+\varepsilon a^{\ast })=f(a)+\varepsilon a^{\ast }f^{\prime }(a)$%
\newline
\noindent where $f^{\prime }(a)$ is the derivation of $f$. Thus we can obtain%
\newline
\indent$sin(a+\varepsilon a^{\ast })=sina+\varepsilon a^{\ast }cosa$ \newline
\indent$cos(a+\varepsilon a^{\ast })=cosa-\varepsilon a^{\ast }sina$\newline
\indent The set of $D^{3}=\{\overrightarrow{A}|\ \ \overrightarrow{A}=%
\overrightarrow{a}+\varepsilon \overrightarrow{a^{\ast }},\overrightarrow{a},%
\overrightarrow{a^{\ast }}\in 
\mathbb{R}
^{3}\}$ is a module on the ring $D$. For any $\overrightarrow{A}=%
\overrightarrow{a}+\varepsilon \overrightarrow{a^{\ast }},\overrightarrow{B}=%
\overrightarrow{b}+\varepsilon \overrightarrow{b^{\ast }}\in D^{3}$, the
scalar or inner product and the vector product of $\overrightarrow{A}$ and $%
\overrightarrow{B}$ are defined by, respectively,

\indent$\langle \overrightarrow{A},\overrightarrow{B}\rangle =\langle 
\overrightarrow{a},\overrightarrow{b}\rangle +\varepsilon (\langle 
\overrightarrow{a},\overrightarrow{b^{\ast }}\rangle +\langle 
\overrightarrow{a^{\ast }},\overrightarrow{b}\rangle )$,\newline
\indent$\overrightarrow{A}\wedge \overrightarrow{B}=\overrightarrow{a}\wedge 
\overrightarrow{b}+\varepsilon (\overrightarrow{a}\wedge \overrightarrow{%
b^{\ast }}+\overrightarrow{a^{\ast }}\wedge \overrightarrow{b}).$\newline
\indent If $\overrightarrow{a}\neq 0$, the norm $\Vert \overrightarrow{A}%
\Vert $ of $\overrightarrow{A}=\overrightarrow{a}+\varepsilon 
\overrightarrow{a^{\ast }}$ is defined by \newline
\indent$\left\Vert \overrightarrow{A}\right\Vert =\sqrt{\left\vert
\left\langle \overrightarrow{A},\overrightarrow{A}\right\rangle \right\vert }%
=\left\Vert \overrightarrow{a}\right\Vert +\varepsilon \frac{\left\langle 
\overrightarrow{a},\overrightarrow{a}^{\ast }\right\rangle }{\left\Vert 
\overrightarrow{a}\right\Vert }\mathrm{\;},\,\mathrm{\;}\left\Vert 
\overrightarrow{a}\right\Vert \neq 0.$ \newline
\indent A dual vector $\overrightarrow{A}$ with norm $1$ is called a dual
unit vector. The set \newline
\indent$S^{2}=\{\overrightarrow{A}=\overrightarrow{a}+\varepsilon 
\overrightarrow{a^{\ast }}\in D^{3}|\Vert \overrightarrow{A}\Vert =(1,0),%
\overrightarrow{a},\overrightarrow{a^{\ast }}\in 
\mathbb{R}
^{3}\}$ \newline
is called the dual unit sphere with the center $\overrightarrow{O}$ in $%
D^{3} $.\newline
\indent Let $\alpha (t)=(\alpha _{1}(t),\alpha _{2}(t),\alpha _{3}(t))$ and $%
\beta (t)=(\beta _{1}(t),\beta _{2}(t),\beta _{3}(t))$ be real valued curves
in $E^{3}$. Then $\widetilde{\alpha }(t)=\alpha (t)+\varepsilon \alpha
^{\ast }(t)$ is a curve in $D^{3}$ and it is called dual space curve. If the
real valued functions $\alpha _{i}(t)$ and $\alpha _{i}^{\ast }(t)$ are
differentiable then the dual space curve $\widetilde{\alpha }(t)$ is
differentiable in $D^{3}$. The real part $\alpha (t)$ of the dual space
curve $\widetilde{\alpha }=\widetilde{\alpha }(t)$ is called indicatrix. The
dual arc-length of real dual space curve $\widetilde{\alpha }(t)$ from $%
t_{1} $ to $t$ is defined by \newline
\indent$\widetilde{s}=\int_{t_{1}}^{t}\Vert \overrightarrow{\widetilde{%
\alpha ^{\prime }}}(t)\Vert dt=\int_{t_{1}}^{t}\Vert \overrightarrow{\alpha {%
^{\prime }}}(t)\Vert dt+\varepsilon ={\int_{t_{1}}^{t}\langle 
\overrightarrow{t},(\overrightarrow{\alpha ^{\ast }}(t))^{^{\prime }}\rangle 
}dt=s+\varepsilon s^{\ast }$

\noindent $\overrightarrow{t}$ is unit tangent vector of the indicatrix $%
\alpha (t)$ which is a real space curve in $IE^{3}$. From now on we will
take the arc length $s$ of $\overrightarrow{\alpha (t)}$ as the parameter
instead of $t$ \newline
\indent The Lorentzian inner product of dual vectors $\overrightarrow{A},%
\overrightarrow{B}\in D^{3}$ is defined by \newline
\indent$\langle \overrightarrow{A},\overrightarrow{B}\rangle =\langle 
\overrightarrow{a},\overrightarrow{b}\rangle +\varepsilon (\langle 
\overrightarrow{a},\overrightarrow{b^{\ast }}\rangle +\langle 
\overrightarrow{a^{\ast }},\overrightarrow{b}\rangle )$ \newline
with the Lorentzian inner product $\overrightarrow{a}=(a_{1},a_{2},a_{3})$
and $\overrightarrow{b}=(b_{1},b_{2},b_{3})\in 
\mathbb{R}
^{3}$ \newline
\indent$\langle \overrightarrow{a},\overrightarrow{b}\rangle
=-a_{1}b_{1}+a_{2}b_{2}+a_{3}b_{3}.$ \newline
Thus, $D^{3},\langle ,\rangle $ is called the dual Lorentzian space and
denoted by $D^{3}$. We call the elements of $D^{3}$ as the dual vectors. For 
$\overrightarrow{A}\neq \overrightarrow{0}$. the norm $\Vert \overrightarrow{%
A}\Vert $ of $\overrightarrow{A}$ is defined by $\left\Vert \overrightarrow{A%
}\right\Vert =\sqrt{\left\vert \left\langle \overrightarrow{A},%
\overrightarrow{A}\right\rangle \right\vert }$ . The dual vector $%
\overrightarrow{A}=\overrightarrow{a}+\varepsilon \overrightarrow{a^{\ast }}$
is called dual spacelike vector if $\left\langle \overrightarrow{A},%
\overrightarrow{A}\right\rangle >0$ or $\overrightarrow{A}=0$, dual timelike
vector if $\left\langle \overrightarrow{A},\overrightarrow{A}\right\rangle
<0 $ , dual lightlike vector if $\left\langle \overrightarrow{A},%
\overrightarrow{A}\right\rangle =0$ for $\overrightarrow{A}\neq 0$. The dual
Lorentzian cross-product of $\overrightarrow{A}\,,\,\overrightarrow{\mathrm{%
\;}B}\in D^{3}\,$is defined by \newline
\indent$\overrightarrow{A}\wedge \overrightarrow{B}=\overrightarrow{a}\wedge 
\overrightarrow{b}+\varepsilon \left( \overrightarrow{a}\wedge 
\overrightarrow{b}^{\ast }+\overrightarrow{a}^{\ast }\wedge \overrightarrow{b%
}\right) $ \newline
where $\overrightarrow{a}\wedge \overrightarrow{b}=\left(
a_{3}b_{2}-a_{2}b_{3},a_{1}b_{3}-a_{3}b_{1},a_{1}b_{2}-a_{2}b_{1}\right) $ $%
\overrightarrow{a},\overrightarrow{b}\in 
\mathbb{R}
^{3}$i s the Lorentzian cross product.

Dual number \textbf{$\Phi =\theta +\varepsilon \theta ^{\ast }$ } is called
dual angle between \textbf{\ $\overrightarrow{A}\,\,\mathrm{ve}\,\,%
\overrightarrow{B}\,\,$}unit dual vectors. Then we was \newline
\indent$\sinh \left( \theta +\varepsilon \theta ^{\ast }\right) =\sinh
\theta +\varepsilon \theta ^{\ast }\cosh \theta $\newline
\indent$\cosh \left( \theta +\varepsilon \theta ^{\ast }\right) =\cosh
\varphi +\varepsilon \theta ^{\ast }\sinh \theta .$ \newline
\indent Let $\left\{ T\left( s\right) ,N\left( s\right) ,B\left( s\right)
\right\} $be the moving Frenet frame along the curve $\widetilde{\alpha }%
\left( s\right) $. Then $T\left( s\right) $,$N\left( s\right) $ and $B\left(
s\right) $ are dual tangent, the dual principal normal and the dual binormal
vector of the curve $\widetilde{\alpha }\left( s\right) $, respectively.
Depending on the casual character of the curve $\widetilde{\alpha }$, we
have the following dual Frenet formulas:

If $\widetilde{\alpha }$ is a dual timelike curve ;

\begin{equation}  \label{GrindEQ__1_1_}
\left(%
\begin{array}{c}
{T^{\prime }} \\ 
{N^{\prime }} \\ 
{B^{\prime }}%
\end{array}%
\right)=\left(%
\begin{array}{ccc}
{0} & {\kappa } & {0} \\ 
{\kappa } & {0} & {\tau } \\ 
{0} & {-\tau } & {0}%
\end{array}%
\right)\left(%
\begin{array}{c}
{T} \\ 
{N} \\ 
{B}%
\end{array}%
\right)
\end{equation}
\newline
where $\left\langle T,T\right\rangle =-1,\left\langle N,N\right\rangle
=\left\langle B,B\right\rangle =1,\, \left\langle T,N\right\rangle
=\left\langle N,B\right\rangle =\left\langle T,B\right\rangle =0.$

\noindent We denote by $\left\{V_{1} \left(s\right),V_{2}
\left(s\right),V_{3} \left(s\right)\right\}$ the moving Frenet frame along
the curve $\widetilde{\beta }\left(s\right)$. Then $V_{1}
\left(s\right),V_{2} \left(s\right)$ and $V_{3} \left(s\right)$ are dual
tangent, the dual principal normal and the dual binormal vector of the curve 
$\widetilde{\beta }\left(s\right)$, respectively. Depending on the casual
character of the curve $\widetilde{\beta }$, we have the following dual
Frenet -- Serret formulas:

If $\widetilde{\beta }$ is a dual timelike curve;

\begin{equation}
\left( 
\begin{array}{c}
{V_{1}^{{^{\prime }}}} \\ 
{V_{2}^{\prime }} \\ 
{V_{3}^{{^{\prime }}}}%
\end{array}%
\right) =\left( 
\begin{array}{ccc}
{0} & {P} & {0} \\ 
{P} & {0} & {Q} \\ 
{0} & -{Q} & {0}%
\end{array}%
\right) \left( 
\begin{array}{c}
{V_{1}} \\ 
{V_{2}} \\ 
{V_{3}}%
\end{array}%
\right)  \label{GrindEQ__1_2_}
\end{equation}%
\newline
where $\left\langle V_{1},V_{1}\right\rangle =-1,\left\langle
V_{2},V_{2}\right\rangle =1,\,\left\langle V_{3},V_{3}\right\rangle
=1,\,\left\langle V_{1},V_{2}\right\rangle =\left\langle
V_{2},V_{3}\right\rangle =\left\langle V_{1},V_{3}\right\rangle =0.$

\noindent If the curves are unit speed curve, then curvature and torsion
calculated by,

\begin{equation}  \label{GrindEQ__1_3_}
\left\{%
\begin{array}{l}
{\kappa =\left\| T^{\prime }\right\| ,} \\ 
{\tau =\left\langle N^{{^{\prime }} } ,B\right\rangle ,} \\ 
{} \\ 
{P=\left\| V_{1} ^{{^{\prime }} } \right\| ,} \\ 
{Q=\left\langle V_{2} ^{{^{\prime }} } ,V_{3} \right\rangle .}%
\end{array}%
\right.
\end{equation}

\noindent If the curves are not unit speed curve, then curvature and torsion
calculated by,

\begin{equation}  \label{GrindEQ__1_4_}
\left\{%
\begin{array}{l}
{\kappa =\frac{\left\| \widetilde{\alpha }^{{^{\prime }} } {\wedge 
\widetilde{\alpha }^{^{\prime \prime }}}\right\|}{{\left\| \widetilde{\alpha 
}^{^{\prime }}{} \right\|}^{3}} ,\, \, \, \, \, \, \, \, \, \, \, \, \, \,
\, \, \, \, \, \, \, \, \, \, \, \, \, \, \, \, \, \, \, \, \, \tau =\frac{%
\det \left(\widetilde{\alpha }^{{^{\prime }} } ,\widetilde{\alpha }^{{%
^{\prime }} {^{\prime }} } ,\widetilde{\alpha }^{{^{\prime }} {^{\prime }} {%
^{\prime }} } \right)}{{\left\| \widetilde{\alpha }^{{^{\prime }} } \wedge 
\widetilde{\alpha }^{{^{\prime }} {^{\prime }} } \right\|}^{2}} ,} \\ 
{} \\ 
{P=\frac{\|{\widetilde{\beta}^{^{\prime }}\wedge{\widetilde{\beta}^{^{\prime
\prime }}\|}}}{\|{\widetilde{\beta}^{^{\prime }}}\|^{3}}},\, \, \, \, \, \,
\, \, \, \, \, \, \, \, \, \, \, \, \, \, \, \, \, \, \, \, \, \, \, \, \,
\, \, \, \, Q =\frac{\det \left(\widetilde{\beta }^{{^{\prime }} } ,%
\widetilde{\beta }^{{^{\prime }} {^{\prime }} } ,\widetilde{\beta }^{{%
^{\prime }} {^{\prime }} {^{\prime }} } \right)}{{\left\| \widetilde{\beta }%
^{{^{\prime }} } \wedge \widetilde{\alpha }^{{^{\prime }} {^{\prime }} }
\right\|}^{2}}%
\end{array}%
\right.
\end{equation}

\noindent \textbf{Definition 2.1. a) Dual Hyperbolic angle: }Let $%
\overrightarrow{A}$ and $\overrightarrow{B}$ be dual timelike vectors in $%
D_{1}^{3}$. Then the dual angle between $\overrightarrow{A}$ and $%
\overrightarrow{B}$ is defined by $\left\langle \overrightarrow{A},%
\overrightarrow{B}\right\rangle =-\left\Vert \overrightarrow{A}\right\Vert
\,\,\left\Vert \overrightarrow{B}\right\Vert \cosh \Phi $. The dual number $%
\Phi =\theta +\varepsilon \theta ^{\ast }$ is called the dual hyberbolic
angle.

\noindent \textbf{b) Dual Central angle: }Let $\overrightarrow{A}$ and $%
\overrightarrow{B}$ be spacelike vectors in $D_{1}^{3}$ that span a dual
timelike vector subspace. Then the dual angle between $\overrightarrow{A}$
and $\overrightarrow{B}$ is defined by $\left\langle \overrightarrow{A},%
\overrightarrow{B}\right\rangle =\left\Vert \overrightarrow{A}\right\Vert
\,\,\left\Vert \overrightarrow{B}\right\Vert \cosh \Phi $. The dual number $%
\Phi =\theta +\varepsilon \theta ^{\ast }$ is called the dual central angle.

\noindent \textbf{c) Dual Spacelike angle: }Let $\overrightarrow{A}$ and $%
\overrightarrow{B}$ be dual spacelike vectors in$D_{1}^{3}$ that span a dual
spacelike vector subspace. Then the dual angle between $\overrightarrow{A}$
and $\overrightarrow{B}$ is defined by $\left\langle \overrightarrow{A},%
\overrightarrow{B}\right\rangle =\left\Vert \overrightarrow{A}\right\Vert
\,\,\left\Vert \overrightarrow{B}\right\Vert \cos \Phi $. The dual number $%
\Phi =\theta +\varepsilon \theta ^{\ast }$ is called the dual spacelike
angle.

\noindent \textbf{d) Dual Lorentzian timelike angle: }Let\textbf{\ }$%
\overrightarrow{A}$ be a dual spacelike vector and $\overrightarrow{B}$ be a
dual timelike vector in $ID_{1}^{3}$. Then the dual angle between $%
\overrightarrow{A}$ and $\overrightarrow{B}$ is defined by $\left\langle 
\overrightarrow{A},\overrightarrow{B}\right\rangle =\left\Vert 
\overrightarrow{A}\right\Vert \,\,\left\Vert \overrightarrow{B}\right\Vert
\sinh \Phi $. The dual number $\Phi =\theta +\varepsilon \theta ^{\ast }$ is
called the dual Lorentzian timelike angle [18, 19, 20].

\section{DUAL TIMELIKE MANNHEIM \newline
PARTNER CURVE IN $D_{1}^{3}$}

\indent In this section, we define dual timelike Mannheim partner curves in $%
D_{1}^{3}$ and we give some characterization for dual timelike Mannheim
partner curves in the same space. Using these relationships, we will comment
again Shell's and Mannheim's theorems.

\noindent \textbf{Definition 3.1. }Let\textbf{\ $\widetilde{\alpha }%
:I\rightarrow ID_{1}^{3}$, $\widetilde{\alpha }\left( s\right) =\alpha
\left( s\right) +\varepsilon \alpha ^{\ast }\left( s\right) $ }and \newline
\textbf{\ $\widetilde{\beta }:I\rightarrow ID_{1}^{3}$,$\widetilde{\beta }%
\left( s\right) =\beta \left( s\right) +\varepsilon \beta ^{\ast }\left(
s\right) $ }be\textbf{\ }dual timelike curves. If there exists a
corresponding relationship between the dual timelike curve $\widetilde{%
\alpha }$ and the dual timelike curve $\widetilde{\beta }$ such that, at the
corresponding points of the curves, the dual binormal lines of $\widetilde{%
\alpha }$ coincides with the dual principal normal lines of $\widetilde{%
\beta }$, then $\widetilde{\alpha }$ is called a dual timelike Mannheim
curve, and $\widetilde{\beta }$ is called a dual Mannheim partner curve of $%
\widetilde{\alpha }$. The pair $\left\{ \widetilde{\alpha },\widetilde{\beta 
}\right\} $ is said to be dual timelike Mannheim pair. Let $\left\{
T,N,B\right\} $ be the dual Frenet frame field along $\widetilde{\alpha }=%
\widetilde{\alpha }\left( s\right) $ and let $\left\{
V_{1},V_{2},V_{3}\right\} $ be the Frenet frame field along $\widetilde{%
\beta }=\widetilde{\beta }\left( s\right) $. On the way $\Phi =\theta
+\varepsilon \theta ^{\ast }$ is dual angle between $T$ and $V_{1}$ , there
is an following equations between the Frenet vectors and their derivative;

\begin{equation}
\left( 
\begin{array}{c}
{V_{1}^{{^{\prime }}}} \\ 
{V_{2}^{{^{\prime }}}} \\ 
{V_{3}^{{^{\prime }}}}%
\end{array}%
\right) =\left( 
\begin{array}{ccc}
\cosh {\Phi } & {\sinh \Phi } & {0} \\ 
{0} & {0} & {1} \\ 
{\sinh \Phi } & {\cosh \Phi } & {0}%
\end{array}%
\right) \left( 
\begin{array}{c}
{T} \\ 
{N} \\ 
{B}%
\end{array}%
\right) .  \label{GrindEQ__2_1_}
\end{equation}%
\newline
\textbf{Theorem 3.1. }The distance between corresponding dual points of the
dual timelike Mannheim partner curves in\textbf{\ }$D_{1}^{3}$ is constant.

\noindent \textbf{Proof: }From the definition of dual spacelike Mannheim
curve, we can write

\begin{equation}  \label{GrindEQ__2_2_}
\tilde{\beta }(s^{*} )=\; \; \tilde{\alpha }(s)+\lambda
\left(s\right)B\left(s\right)
\end{equation}

\noindent By taking the derivate of this equation with respect to $s$ and
applying the Frenet formulas, we get

\begin{equation}  \label{GrindEQ__2_3_}
V_{1} \frac{ds^{*} }{ds} =T-\lambda \tau N+\lambda ^{\prime }B
\end{equation}
where the superscript $\left(^{\prime }\right)$ denotes the derivative with
respect to the arc length parameter s of the dual curve $\tilde{\alpha }(s)$%
. Since the dual vectors $B$ and $V_{2} $ are linearly, we get

\indent $\left\langle V_{1} \frac{ds^{*} }{ds} ,B\right\rangle =\left\langle
T,B\right\rangle -\lambda \tau \left\langle N,B\right\rangle +\lambda
^{\prime }\left\langle B,B\right\rangle$ and $\lambda ^{\prime }=0$ \newline
If we take $\lambda =\lambda _{1} +\varepsilon \lambda _{1}^{*} $, we get $%
\lambda ^{\prime }_{1} =0$ ve $\lambda _{1}^{*^{\prime }} =0$ . From here,
we can write $\lambda _{1} =c_{1}$ and $\lambda_{1}^{*}=c_{2},$ $%
c_{1},c_{2}=cons.$

\noindent Then we get $\lambda =c_{1} +\varepsilon c_{2} $. On the other
hand, from the definition\textbf{\ }of distance function between $\tilde{%
\alpha }(s)$ and $\tilde{\beta }(s)$ we can write \newline
\indent $d\left(\tilde{\alpha }(s),\tilde{\beta }(s)\right)=\left\| \tilde{%
\beta }(s)-\tilde{\alpha }(s)\right\| = \left|\lambda _{1} \right|\mp
\varepsilon \lambda _{1}^{*} = \left|c_{1} \right|\mp \varepsilon c_{2}$ 
\newline
This is completed the proof.

\noindent \textbf{Theorem 3.2. }For a dual timelike curve\textbf{\ }$%
\widetilde{\alpha }$ in \textbf{\ }$D_{1}^{3}$, there is a dual timelike
curve $\tilde{\beta}$ so that $\left\{ \tilde{\alpha},\tilde{\beta}\right\} $
is a dual timelike Mannheim pair.\newline
\noindent \textbf{Proof: }Since\textbf{\ }the dual vectors $V_{2}$ and $B$
are linearly dependent, the equation (3.2) can be written as

\begin{equation}  \label{GrindEQ__2_4_}
\tilde{\alpha }=\tilde{\beta }-\lambda V_{2}
\end{equation}
Since $\lambda $ is a nonzero constant, there is a dual timelike curve $%
\tilde{\beta }$ for all values of $\lambda $.

Now, we can give the following theorem related to curvature and torsion of
the dual timelike Mannheim partner curves.

\noindent \textbf{Theorem 3.3. }Let $\left\{ \tilde{\alpha},\tilde{\beta}%
\right\} $ be a dual timelike Mannheim pair in $D_{1}^{3}$\textbf{. }If%
\textbf{\ }$\tau $ is dual torsion of \textbf{\ }$\tilde{\alpha}$ and $P$ is
dual curvature and $Q$ is dual torsion of $\tilde{\beta}$ , then

\begin{equation}
\tau =\frac{P}{\lambda Q}  \label{GrindEQ__2_5_}
\end{equation}%
\textbf{Proof: }By taking the derivate of equation\textbf{\ }(3.3) with
respect to $s$ and applying the Frenet formulas, we obtain

\begin{equation}
V_{1}\frac{ds^{\ast }}{ds}=T-\lambda \tau N  \label{GrindEQ__2_6_}
\end{equation}%
Let \textbf{$\Phi =\theta +\varepsilon $}$\theta ^{\ast }$\textbf{\ }be%
\textbf{\ }dual angle between\textbf{\ } the dual tangent vectors $T$ and $%
V_{1}$, we can write

\begin{equation}
\left\{ 
\begin{array}{l}
{V_{1}=\cosh \Phi \,T+\sinh \Phi \,N} \\ 
{V_{3}=\sinh \Phi \,T+\cosh \Phi \,N}%
\end{array}%
\right.  \label{GrindEQ__2_7_}
\end{equation}%
From (3.6) and (3.7) , we get

\begin{equation}
\frac{ds^{\ast }}{ds}=\frac{1}{\cosh \Phi },\,\,\,\,\,-\lambda \tau =\sinh
\Phi \frac{ds^{\ast }}{ds}  \label{GrindEQ__2_8_}
\end{equation}%
By taking the derivate of equation\textbf{\ }(3.4) with respect to $s$ and
applying the Frenet formulas, we obtain

\begin{equation}
T=\left( 1-\lambda P\right) V_{1}\frac{ds^{\ast }}{ds}-\lambda QV_{3}\frac{%
ds^{\ast }}{ds}  \label{GrindEQ__2_9_}
\end{equation}%
From equation (3.7) we can write

\begin{equation}
\left\{ 
\begin{array}{l}
{T=\cosh \,V_{1}-\sinh \Phi \,V_{3}} \\ 
{N=-\sinh \Phi \,V_{1}+\cosh \Phi \,V_{3}}%
\end{array}%
\right.  \label{GrindEQ__2_10_}
\end{equation}%
where $\Phi $ is the dual angle between $T$ and $V_{1}$ at the corresponding
points of the dual curves of $\tilde{\alpha}$ and $\tilde{\beta}$ . By
taking into consideration equations (3.9) and (3.10), we get

\begin{equation}
\cosh \Phi =\left( 1-\lambda P\right) \frac{ds^{\ast }}{ds},\,\,\,\,\,\sinh
\Phi =\lambda Q\frac{ds^{\ast }}{ds}  \label{GrindEQ__2_11_}
\end{equation}%
Substituting $\frac{ds^{\ast }}{ds}$ into (3.11) , we get

\begin{equation}
\cosh ^{2}\Phi =-\left( 1+\lambda P\right) ,\,\,\,\,\,\sinh ^{2}\Phi
=\lambda ^{2}\tau Q  \label{GrindEQ__2_12_}
\end{equation}%
From the last equation, we can write \newline
\indent$\tau =\frac{P}{\lambda Q}$ \newline
If the last equation is seperated into the dual and real parts, we can obtain

\begin{equation}
\left\{ 
\begin{array}{l}
{k_{2}=\frac{p}{cq}} \\ 
{k_{2}^{\ast }=\frac{p^{\ast }q-pq^{\ast }}{cq^{2}}}%
\end{array}%
\right.  \label{GrindEQ__2_13_}
\end{equation}%
\textbf{Corollary 3.1. }Let $\left\{ \tilde{\alpha},\tilde{\beta}\right\} $
be a dual timelike Mannheim pair in $D_{1}^{3}$\textbf{. }Then, the dual
product of torsions\textbf{\ }$\tau $ and $Q$ at the corresponding points of
the dual timelike Mannheim partner curves is not constant.

Namely, Schell's theorem is invalid for the dual timelike Mannheim curves.
By considering Theorem 3.3 we can give the following results.

\noindent \textbf{Corollary 3.2. }Let $\left\{ \tilde{\alpha},\tilde{\beta}%
\right\} $ be a dual timelike Mannheim pair in $D_{1}^{3}$\textbf{. }Then,%
\textbf{\ }torsions\textbf{\ }$\tau $\textit{\ }and $Q$ has a negative sign.

\noindent \textbf{Theorem 3.4.} Let $\left\{ \tilde{\alpha},\tilde{\beta}%
\right\} $ be a dual timelike Mannheim pair in $D_{1}^{3}$\textbf{. }Between
the curvature and the torsion of the dual timelike curve $\widetilde{\beta }$
, there is the relationship

\begin{equation}
\mu Q+\lambda P=1  \label{GrindEQ__2_14_}
\end{equation}%
where $\mu $ and$\lambda $ are nonzero dual numbers.

\noindent \textbf{Proof: }From equation (3.11), we obtain \newline
\indent$\frac{\cosh \Phi }{1-\lambda P}=\frac{\sinh \Phi }{\lambda Q}$, 
\newline
arranging this equation, we get \newline
\indent$\tanh \Phi =\frac{1-\lambda P}{\lambda Q}$, \newline
and if we choose $\mu =\lambda \tanh \Phi $ for brevity, we see that \newline
\indent$\mu Q+\lambda P=1$.\newline
\textbf{Theorem 3.5. }Let $\left\{ \tilde{\alpha},\tilde{\beta}\right\} $ be
a dual timelike Mannheim pair in $D_{1}^{3}$\textbf{. }There are the
following equations for the curvatures and the torsions of the curves $%
\widetilde{\alpha }$ ve $\widetilde{\beta }$ \newline
\indent$i)\kappa =-\frac{d\Phi }{ds},$ \newline
\indent$ii)\tau =-P\sinh \Phi \frac{ds^{\ast }}{ds}-Q\cosh \Phi \frac{%
ds^{\ast }}{ds},$ \newline
\indent$iii)P=\tau \sinh \Phi \frac{ds}{ds^{\ast }},$ \newline
\indent$iv)Q=-\tau \cosh \Phi \frac{ds}{ds^{\ast }}.$

\textbf{Proof: }$i)$By considering equation (3.7), we can easily that $%
\left\langle T,V_{1}\right\rangle =-\cosh \Phi $. Differentiating of this
equality with respect to \textit{s }by considering equation (2.1) , we have 
\newline
\indent$\left\langle T^{\prime },V_{1}\right\rangle +\left\langle T,V_{1}^{{%
^{\prime }}}\right\rangle =-\sinh \Phi \frac{d\Phi }{ds}$,\newline
from equations (2.1) and (2.2), we can write \newline
\indent$\left\langle \kappa N,V_{1}\right\rangle +\left\langle T,PV_{2}\frac{%
ds^{\ast }}{ds}\right\rangle =-\sinh \Phi \frac{d\Phi }{ds}$, \newline
from equations (3.10), we get \newline
\indent$\kappa =-\frac{d\Phi }{ds}$.\newline
If the last equation is seperated into the dual and real part, we can obtain 
\newline
$ii)$ By considering equation (3.7), we can easily that $\left\langle
N,V_{2}\right\rangle =0$. Differentiating of this equality with respect to 
\textit{s and }by considering equation (2.1) , we have \newline
\indent$\left\langle N^{\prime },V_{2}\right\rangle +\left\langle N,V_{2}^{{%
^{\prime }}}\frac{ds^{\ast }}{ds}\right\rangle =0$, \newline
From equations (2.1) and (2.2), we can write \newline
\indent$\left\langle \kappa T+\tau B,V_{2}\right\rangle +\left\langle -\sinh
\Phi \,V_{1}+\cosh \Phi \,V_{3},\left( PV_{1}+QV_{3}\right) \frac{ds^{\ast }%
}{ds}\right\rangle =0$,\newline
From equations (3.10), we get \newline
\indent$\tau =-P\sinh \Phi \frac{ds^{\ast }}{ds}-Q\cosh \Phi \frac{ds^{\ast }%
}{ds}$, \newline
$iii)$ By considering equation (3.7), we can easily that $\,\left\langle
B,V_{1}\right\rangle =0$. Differentiating of this equality with respect to 
\textit{s and }by considering equation (2.1), we have \newline
\indent$\left\langle B^{\prime },V_{1}\right\rangle +\left\langle B,V_{1}^{{%
^{\prime }}}\frac{ds^{\ast }}{ds}\right\rangle =0$, \newline
From equations (2.1), (2.2) and (3.10) we can write \newline
\indent$\left\langle -\tau \left( -\sinh \Phi \,V_{1}+\cosh \Phi
\,V_{3}\right) ,V_{1}\right\rangle +\left\langle B,PV_{2}\frac{ds^{\ast }}{ds%
}\right\rangle =0$,\newline
\indent$P=\tau \sinh \Phi \frac{ds}{ds^{\ast }}$,\newline
$iv)$ By considering equation (3.7), we can easily that $\left\langle
B,V_{3}\right\rangle =0$. Differentiating of this equality with respect to 
\textit{s }by considering equation (2.1) , we have \newline
\indent$\left\langle B^{\prime },V_{3}\right\rangle +\left\langle B,V_{3}^{{%
^{\prime }}}\frac{ds^{\ast }}{ds}\right\rangle =0$,\newline
From equations (2.1), (2.2) and (3.10) we can write \newline
\indent$\left\langle -\tau \left( -\sinh \Phi \,V_{1}+\cosh \Phi
\,V_{3}\right) ,V_{3}\right\rangle +\left\langle B,-QV_{2}\frac{ds^{\ast }}{%
ds}\right\rangle =0$,\newline
\indent$Q=-\tau \cosh \Phi \frac{ds}{ds^{\ast }}$. \newline
\textbf{Corollary 3.3. }Let $\left\{ \tilde{\alpha},\tilde{\beta}\right\} $
be a dual timelike - spacelike Mannheim pair in $D_{1}^{3}$\textbf{. }If the
statements of Theorem 3.5 is seperated into the dual and real part, we can
obtain \newline
\indent$i){\ \left\{ 
\begin{array}{c}
{k_{2}=-p\sinh \theta \frac{ds^{\ast }}{ds}-q\cosh \theta \frac{ds^{\ast }}{%
ds}}\newline
\\ 
{k_{2}^{\ast }=-\left( p^{\ast }\sinh \theta +p\theta ^{\ast }\cosh \theta
\right) \frac{ds^{\ast }}{ds}-\left( q^{\ast }\cosh \theta +q\theta ^{\ast
}\sinh \theta \right) \frac{ds^{\ast }}{ds}}%
\end{array}%
\right. }$

\begin{equation*}
\,ii)\left\{ 
\begin{array}{l}
{p=k_{2}\sinh \theta \frac{ds}{ds^{\ast }}} \\ 
{p^{\ast }=\left( k_{2}^{\ast }\sinh \theta +k_{2}\theta ^{\ast }\cosh
\theta \right) \frac{ds}{ds^{\ast }},}%
\end{array}%
\right.
\end{equation*}

\begin{equation*}
iii)\left\{ 
\begin{array}{l}
{q=-k_{2}\cosh \theta \frac{ds}{ds^{\ast }}} \\ 
{q^{\ast }=-\left( k_{2}^{\ast }\cosh \theta +k_{2}\theta ^{\ast }\theta
\right) \frac{ds}{ds^{\ast }}.}%
\end{array}%
\right.
\end{equation*}%
By considering the statements iii and iv) of Theorem 2.5 we can give the
following results.

\noindent \textbf{Corollary 3.4. }Let $\left\{ \tilde{\alpha},\tilde{\beta}%
\right\} $ be a dual timelike Mannheim pair in $D_{1}^{3}$\textbf{. }Then
there exist the following relation between curvature and torsion of\textbf{\ 
}$\widetilde{\beta }$ and torsion of $\widetilde{\alpha }$;

\begin{equation}
Q^{2}-P^{2}=\tau ^{2}\left( \frac{ds}{ds^{\ast }}\right) ^{2}
\label{GrindEQ__2_15_}
\end{equation}%
\textbf{Theorem 3.6. }A dual timelike space curve in\textbf{\ }$ID_{1}^{3}$
is a dual timelike Mannheim curve if and only if its curvature $P$ and
torsion $Q$ satisfy the formula

\begin{equation}
\lambda \left( P^{2}-Q^{2}\right) =P  \label{GrindEQ__2_16_}
\end{equation}%
where $\lambda $ is never pure dual constant.

\noindent \textbf{Proof:} By taking the derivate of the statement $%
\widetilde{\alpha }=\widetilde{\beta }-\lambda V_{2}$ with respect to $s$
and applying the Frenet formulas we obtain \newline
\indent$T\frac{ds}{ds^{\ast }}=V_{1}-\lambda \left( PV_{1}+QV_{3}\right) $, 
\newline
\indent$\kappa N\left( \frac{ds}{ds^{\ast 2}}\right) +T\frac{d^{2}s}{%
ds^{\ast 2}}=PV_{2}-\lambda \left( P^{\prime }V_{1}+Q^{\prime }V_{3}+\left(
P^{2}-Q^{2}\right) V_{2}\right) $ \newline
Taking the inner product the last equation with $B$, we get \newline
\indent$\lambda \left( P^{2}-Q^{2}\right) =P$. \newline
If the last equation is seperated into the dual and real part, we can obtain

\begin{equation}
\left\{ 
\begin{array}{l}
{p=\lambda \left( p^{2}-q^{2}\right) } \\ 
{p^{\ast }=2\lambda \left( pp^{\ast }-qq^{\ast }\right) }%
\end{array}%
\right.  \label{GrindEQ__2_17_}
\end{equation}%
where $\lambda =c_{1}+\varepsilon c_{2}$ .

\noindent \textbf{Theorem 3.7. }Let $\left\{ \tilde{\alpha},\tilde{\beta}%
\right\} $ be a dual timelike Mannheim partner curves in $D_{1}^{3}$\textbf{%
. }Moreover, the dual points\textbf{\ }$\widetilde{\alpha }\left( s\right) $%
, $\widetilde{\beta }\left( s\right) $ be two corresponding dual points of $%
\left\{ \widetilde{\alpha },\widetilde{\beta }\right\} $and $M$ ve $M^{\ast
} $ be the curvature centers at these points, respectively. Then, the ratio 
\newline

\begin{equation}  \label{GrindEQ__2_18_}
\frac{\left\| \widetilde{\beta }\left(s\right)M\right\| }{\left\| \widetilde{%
\alpha }\left(s\right)M\right\| } :\frac{\left\| \widetilde{\beta }%
\left(s\right)M^{*} \right\| }{\left\| \widetilde{\alpha }%
\left(s\right)M^{*} \right\| }=\left(1+\kappa P\right)\left(1+\lambda
P\right)\ne constant.
\end{equation}

\noindent \textbf{Proof: }A circle that lies in the dual osculating plane of
the point $\widetilde{\alpha }\left( s\right) $ on the dual timelike curve $%
\widetilde{\alpha }$ and that has the centre $M=\widetilde{\alpha }\left(
s\right) +\frac{1}{\kappa }N$ lying on the dual principal normal $N$ of the
point $\widetilde{\alpha }\left( s\right) $ and the radius $\frac{1}{\kappa }
$ far from $\widetilde{\alpha }\left( s\right) $, is called dual osculating
circle of the dual curve $\widetilde{\alpha }$ in the point $\widetilde{%
\alpha }\left( s\right) $. Similar definition can be given for the dual
curve $\widetilde{\beta }$ too.

\noindent Then, we can write \newline
\indent $\left\| \widetilde{\alpha }\left(s\right)M\right\| =\left\| \frac{1%
}{\kappa } N\right\| \, \, =\frac{1}{\kappa }$,\newline
\indent $\left\| \widetilde{\alpha }\left(s\right)M^{*} \right\| =\left\|
\lambda B+\frac{1}{P} V_{2} \right\| =\frac{1}{P} +\lambda$, \newline
\indent $\left\| \widetilde{\beta }\left(s\right)M^{*} \right\| =\left\| 
\frac{1}{P} V_{2} \right\| =\frac{1}{P}$, \newline
\indent $\left\| \widetilde{\beta }\left(s\right)M\right\| =\left\| \lambda
V_{3} +\frac{1}{\kappa } N\right\| =\frac{1}{\kappa } +\lambda $ \newline
Therefore, we obtain \newline
\indent $\frac{\left\| \widetilde{\beta }\left(s\right)M\right\| }{\left\| 
\widetilde{\alpha }\left(s\right)M\right\| } :\frac{\left\| \widetilde{\beta 
}\left(s\right)M^{*} \right\| }{\left\| \widetilde{\alpha }%
\left(s\right)M^{*} \right\| } =\left(1+\lambda P\right)\sqrt{1-\lambda ^{2}
\kappa ^{2} } \ne cons. $ \newline
Thus, we can give the following

\noindent \textbf{Corollary 3.5.} Mannheim's Theorem is invalid for the dual
timelike Mannheim partner curve $\left\{ \tilde{\alpha},\tilde{\beta}%
\right\} $in $D_{1}^{3}$\textbf{.}

\end{document}